\documentclass[12pt]{amsart}

\newtheorem{theorem}{Theorem}[section]
\newtheorem{lemma}[theorem]{Lemma}
\newtheorem{corollary}[theorem]{Corollary}
\newtheorem{proposition}[theorem]{Proposition}

\newtheorem{Question}[theorem]{Question}
\newcommand{\ncom}{\newcommand}
\ncom{\rar}{\rightarrow}
\ncom{\lrar}{\longrightarrow}
\ncom{\ov}{\overline}
\ncom{\m}{\mbox}
\ncom{\sta}{\stackrel}
\ncom{\comx}{{\mathbb C}}
\ncom{\Z}{{\mathbb Z}}
\ncom{\Q}{{\mathbb Q}}
\ncom{\R}{{\mathbb R}}
\ncom{\G}{{\mathbb G}}
\ncom{\al}{\alpha}
\ncom{\p}{{\mathbb P}}
\ncom{\E}{{\mathbb E}}
\ncom{\N}{{\mathbb N}}
\ncom{\K}{{\mathbb K}}
\ncom{\Le}{{\mathbb L}}
\ncom{\A}{{\mathbb A}}
\ncom{\B}{{\mathbb B}}
\ncom{\F}{{\mathbb F}}
\ncom{\C}{{\mathbb C}}
\ncom{\f}{\frac}
\ncom{\cA}{{\mathcal A}}
\ncom{\cX}{{\mathcal X}}
\ncom{\cO}{{\mathcal O}}
\ncom{\cW}{{\mathcal W}}
\ncom{\cL}{{\mathcal L}}
\ncom{\cP}{{\mathcal P}}
\ncom{\cH}{{\mathcal H}}
\ncom{\cS}{{\mathcal S}}
\ncom{\cM}{{\mathcal M}}
\ncom{\cC}{{\mathcal C}}
\ncom{\cT}{{\mathcal T}}
\ncom{\cF}{{\mathcal F}}
\ncom{\cN}{{\mathcal N}}
\ncom{\cJ}{{\mathcal J}}
\ncom{\cV}{{\mathcal V}}
\ncom{\cZ}{{\mathcal Z}}
\ncom{\cU}{{\mathcal U}}
\ncom{\cSU}{{\mathcal S \mathcal U}}
\ncom{\cG}{{\mathcal G}}
\ncom{\cQ}{{\mathcal Q}}
\ncom{\cR}{{\mathcal R}}
\ncom{\cE}{{\mathcal E}}
\ncom{\cY}{{\mathcal Y}}
\input xy
\xyoption{all}

\begin{document}
\baselineskip=16pt

\title[Note on unirationality]{A Note on the unirationality of a moduli space of double covers}

\author[J. N. Iyer]{Jaya NN Iyer}
\author[S. M\"uller--Stach]{Stefan M\"uller--Stach}
\address{The Institute of Mathematical Sciences, CIT
Campus, Taramani, Chennai 600113, India}
\address{Department of Mathematics and Statistics, University of Hyderabad, Gachibowli, Central University P O, Hyderabad-500046, India}
\email{jniyer@imsc.res.in}

\address{Mathematisches Institut der Johannes Gutenberg University\"at Mainz,
Staudingerweg 9, 55099 Mainz, Germany}
\email{mueller-stach@uni-mainz.de}

\footnotetext{This work is partly supported by 
Sonderforschungsbereich/Transregio 45.}

\footnotetext{Mathematics Classification Number: 14D05, 14D20, 14C25, 
 14D21}
\footnotetext{Keywords: Moduli spaces, curves, algebraic groups, Chow groups}

\begin{abstract} In this note we look at the moduli space $\cR_{3,2}$ of double covers of genus three curves, branched along $4$ distinct 
points. This space was studied by Bardelli, Ciliberto and Verra in \cite{BCV}. It admits a dominating morphism $\cR_{3,2} \to {\mathcal A}_4$ 
to Siegel space. We show that there is a birational model of $\cR_{3,2}$ as a group quotient of a product of two Grassmannian varieties. 
This gives a proof of the unirationality of $\cR_{3,2}$ and hence a new proof for the unirationality of ${\mathcal A}_4$. 
\end{abstract}
\maketitle


\setcounter{tocdepth}{1}

\section{Introduction}
In this short note we study the moduli space $\cR_{3,2}$ of Bardelli-Ciliberto-Verra \cite{BCV}. 
$\cR_{3,2}$ parametrises the following data: triples $(C,B,L)$, where $C$ is a smooth connected projective curve of genus $3$, 
$L$ is a line bundle of degree $2$ on $C$ and $B$ is a divisor in the linear system $|L^2|$, 
consisting of distinct points. Since a double cover of a genus $3$ curve ramified in $4$ points has genus $7$, this moduli
problem may be viewed as parametrising certain generalized Prym varieties of dimension $4$ which are quotients of a genus 
$7$ Jacobian \cite{BCV}.  The importance of $\cR_{3,2}$ comes from the fact that this Prym construction provides a dominant, 
generically finite morphism $\cR_{3,2} \to {\mathcal A}_4$ to Siegel space, 
which implies in particular that the generic principally polarised abelian fourfold contains 
a finite number of curves of the minimal genus $7$ \cite{BCV}. 

In this note we first describe $\cR_{3,2}$ as follows:

\begin{theorem}\label{uni}
The moduli space $\cR_{3,2}$ is birational to the group quotient 
$$
\cR_{3,2}\sim (G(3,U^+)\times G(4,U^-))/H
$$
of a product of Grassmannians $G(3,U^+)\times G(4,U^-)$ 
by a certain subgroup $H\subset SO(10)$ which is contained in the centraliser of the action of an involution $i$ on $SO(10)$. 
Moreover, there is an irreducible $16$-dimensional projective representation $U$ of $SO(10)$ and $U=U^+\oplus U^-$ is a 
splitting as $\pm$-eigenspaces for the involution $i$ acting on $U$.
\end{theorem}

This description of $\cR_{3,2}$ is similar to the descriptions obtained for the various moduli spaces $\cM_g$, for small $g\leq 9$, 
by Mukai and others (for example, see \cite{Mukai3}, \cite{Mukai6}).

In another direction, it is of wide interest to know when such moduli spaces are rational or unirational varieties. It is known from the 
results of Severi, Sernesi, Katsylo, Mukai, Dolgachev, Chang-Ran, Verra that the moduli spaces $\cM_g$, for small $g\leq 14$ are 
unirational \cite{Sernesi}, \cite{Katsylo}, \cite{Dolgachev}, \cite{Verra}, \cite{chen}. Some moduli spaces of bielliptic covers have also been 
shown to be rational by Bardelli-Del Centina \cite{B-dC}. Unirationality of the moduli space of \'etale double covers of genus $5$ curves, $\cR_5$, is known, due to  Izadi-Lo Giudice-Sankaran   \cite{Izadi} and Verra \cite{Verra2}. See also a very recent preprint
\cite{Fa-Ve}.
 
The above description of $\cR_{3,2}$ in Theorem \ref{uni} says 

\begin{theorem}
The moduli space $\cR_{3,2}$ (and hence ${\mathcal A}_4$) is a unirational variety.
\end{theorem}

Although the implication about ${\mathcal A}_4$ was known before (due to Clemens), our method gives some interesting insight and may have further applications. 
It would be interesting to study the subgroup $H$ in the theorem and its representations $U^+$ and $U^-$ in more detail. We pose some 
open questions in this direction, see \ref{question}.

The proof of the main theorem proceeds by analysing Mukai's description of the moduli space $\cM_7$ and 
restricting the attention to the sublocus $\cR_{3,2}\subset \cM_7$. This sublocus is contained in the singular locus of $\cM_7$ and 
parametrises curves with an involution. The involution plays a crucial role in determining the Grassmannian varieties, in the 
statement of Theorem \ref{uni}.

{\Small Acknowledgements: 
The first named author acknowledges and thanks the Women in 
Mathematics Program on 'Algebraic geometry and Group Actions' in May 2007, at IAS Princeton. 
She also thanks the Maths Department at Mainz, for their 
hospitality and support in June 2008 when this work was partly done. We thank A. Verra for interesting communications 
and B. Totaro for pointing out some errors and for suggestions. }


\section{A birational model of the moduli space $\cR_{3,2}$ of Bardelli-Ciliberto-Verra}\label{birationalmodel}


In this section, we will look at the moduli space
$\cR_{3,2}$ studied by Bardelli-Ciliberto-Verra \cite{BCV}. 
More precisely, let $\cR_{3,2}$ be the moduli space of all isomorphism classes of double coverings $f:C'\rar C$ with $C$ a smooth curve 
of genus $3$, $C'$ irreducible and $f$  branched at $4$ distinct points of $C$. Alternatively, $\cR_{3,2}$ is the moduli space of 
isomorphism classes of triples $(C,B,L)$, where $C$ is a smooth curve of genus $3$, $B$ is an effective divisor on $C$ formed by
$4$ distinct points and $L$ is a line bundle on $C$ such that $L^{\otimes 2}\simeq \cO(B)$.

Note that the genus of the curve $C'$ is $g'=7$ and $\cR_{3,2}\subset \cM_{7}$ (see \cite[p.138]{Cornalba2}).
Then we have
$$
\mbox{dim}\cR_{3,2}= 10.
$$
Our main theorem in this section is the following:

\begin{theorem}\label{birational}

The moduli space $\cR_{3,2}$ is birational to the group quotient of a product of Grassmannians $G(3,U^+)\times G(4,U^-)$, by an 
algebraic subgroup $H\subset SO(10)$.  Here $H$ is contained in the centraliser of the action of an involution $i$ on $SO(10)$. 
Moreover, there is an irreducible $16$-dimensional projective representation $U$ of $SO(10)$ and $U=U^+\oplus U^-$ is a splitting 
as $\pm$-eigenspaces for the involution $i$ acting on $U$.
\end{theorem}

Our proof follows by analysing  Mukai's classification \cite{Mukai3}, \cite{Mukai6} of the generic genus $7$ canonical curve, 
taking into account the action of the involution. Whenever a genus $7$ smooth curve is not tetragonal, then it is a linear 
section of an orthogonal Grassmannian $X_{10}\subset \p^{15}$, given by the spinor embedding (see \cite[p.1632]{Mukai3}). 
Here $\p^{15}=\p(U_{16})$ where $U_{16}$ is the irreducible spinor representation of the spin group $Spin(10)$. 
Hence the space $U_{16}$ is a projective representation of the special orthogonal group $SO(10)$. Projectively, 
this can be translated to say that the group $SO(10)$ acts on $\p^{15}$ and leaves the orthogonal Grassmannian $X_{10}$ 
invariant. In particular $SO(10)$ also acts on the linear subspaces of $\p^{15}$ and we will require its action on the 
Grassmannian $G(7,U_{16})$. This is because a general linear subspace $\p^6\subset \p^{15}$ restricted to $X_{10}$
gives a canonical curve $C$ of genus $7$. In other words, $\p^6$ is the complete linear system given by the canonical bundle 
on $C=\p^6\cap X_{10}$.

Furthermore, we have the following result on the embedding into the homogeneous space.

\begin{theorem}\label{automorphism}
Assume that two linear spaces $P_1,P_2$ cut out smooth curves $C_1, C_2$ from the symmetric space $X_{10}\subset \p^{15}$ respectively.
Then any isomorphism from $C_1$ onto $C_2$ extends to an automorphism $\phi$ of $X_{10}\subset \p^{15}$ with $\phi(P_1)=P_2$.
\end{theorem}
\begin{proof}
See \cite[Theorem 3]{Mukai6}.
\end{proof}
This theorem characterises the non-tetragonal curves of genus $7$.
Explicitly, the moduli space has the following birational model \cite[\S 5, p.1639]{Mukai3}:
\begin{eqnarray*}
\cM_7 & \sim & G(7,U_{16})/SO(10).\\
\end{eqnarray*}

To obtain a birational model of $\cR_{3,2}$, we will utilise the above birational model of $\cM_7$ and analyse the birational equivalence 
restricted to the sublocus $\cR_{3,2}$.

We will need the following lemma in our proof of Theorem \ref{birational}.
We say that a curve $C'$ is tetragonal if and only if there is a line bundle $\cL\in g^1_4(C')$.

Note  that the data $(C,B,L)\in \cR_{3,2}$ corresponds to the data $(C',i)$, where $C'$ is a genus $7$ curve with an involution $i$. Denote
the quotient map $f: C'\rar C=C'/<i>$.

\begin{lemma}\label{nontetragonal} With notations as above, 
consider a double cover $f:C'\rar C$, defined by a line bundle $L$ branched along the set $B$ of $4$ distinct points, and such that 
$L^2=\cO(B)$. Assume that $C, C'$ are not hyperelliptic.
The curve $C'$ has a $\cL\in g^1_4$ only if $\cL$ is the pullback of a line bundle of degree $2$ on $C$.
\end{lemma}
\begin{proof}
The arguments are similar to \cite[Proposition 2.5, p.234]{Ramanan}, and we explain them below.
Let $\cL\in g^1_4(C')$, i.e., $\cL$ is a line bundle of degree $4$ on $C'$ and $h^0(\cL)= 2$.  If $\cL\simeq i^*L$ then $\cL$ descends down to the 
quotient curve $C$ as a line bundle of degree $2$. Suppose $\cL$ is not isomorphic to $i^*\cL$.
Consider the evaluation sequence:
$$
0\rar N \rar H^0(\cL)\otimes \cO_{C'} \rar \cL\rar 0.
$$
Since $h^0(\cL)=2$ we see that $N\simeq \cL^{-1}$.
Tensor the above exact sequence by $i^*\cL$ and take its global sections. Since $\cL\neq i^*\cL$, we observe that $H^0(N\otimes i^*\cL)=0$ and 
hence $H^0(\cL)\otimes H^0(i^*\cL)\subset H^0(\cL\otimes i^*\cL)$. In particular,
$h^0(\cL\otimes i^*\cL)\geq 4$. Since $C'$ is non-hyperelliptic, by Clifford's theorem \cite[IV,5.4]{Arbarello},
$h^0(\cL\otimes i^*\cL)\leq 4$. Hence we obtain the equality $H^0(\cL)\otimes H^0(i^*\cL)\,=\, H^0(\cL\otimes i^*\cL)$.

Now, notice that the line bundle $\cL\otimes i^*\cL$ has degree $8$ on $C'$ and is invariant under $i$. Hence the product line bundle descends 
down to $C$ as a line bundle of degree $4$. Call this line bundle $M$. In other words, $\cL\otimes i^*\cL\simeq f^*M$.
Consider the direct image
$$
f_\ast(\cO_{C'})= \cO_C\oplus L^{-1}.
$$
Hence, by the projection formula, $f_*(\cL\otimes i^*\cL)= M\oplus (M\otimes L^{-1})$.
This gives a decomposition 
$$
H^0(C',\cL\otimes i^*\cL)=H^0(C,M)\oplus H^0(C,M\otimes L^{-1}).
$$
Moreover, we can identify the eigenspaces for the involution $i$ as follows:
\begin{equation}\label{eigen}
H^0(C',\cL\otimes i^*\cL)^+\,=\,H^0(C,M),\,\, H^0(C',\cL\otimes i^*\cL)^-=H^0(C,M\otimes L^{-1}).
\end{equation}
By Riemann-Roch applied to $M$ and $M\otimes L^{-1}$ on $C$, we get the dimension counts: $h^0(M)=3$ if $M=\omega_C$, otherwise $h^0(M)=2$.
Furthermore, since $C$ is non-hyperelliptic 
\begin{equation}\label{clifford}
h^0(M\otimes L^{-1})= 0.
\end{equation} 
by Clifford's theorem and Riemann-Roch.
This implies that
\begin{equation}\label{equidim}
H^0(\cL)\otimes H^0(i^*\cL)=H^0(\cL\otimes i^*\cL)=H^0(f^*M)=H^0(M).
\end{equation}
The first equality in \eqref{equidim} implies that the $\pm$-eigenspaces for the involution $i$ are non-zero. This gives a contradiction 
to \eqref{eigen} and \eqref{clifford}.

\end{proof}

\begin{corollary}\label{generic}
The generic curve in $\cR_{3,2}$ is non-tetragonal.
\end{corollary}
\begin{proof}
By formula~\eqref{clifford} in the proof of Lemma \ref{nontetragonal}, the generic line bundle $\cM$ of degree $2$ on a generic curve of 
genus $3$ has no section. The eigenspace decomposition for the sections of the pullback bundle $\cL':=f^*\cM$ is given as
$$
H^0(C',\cL')=H^0(C,\cM)\oplus H^0(C,\cM\otimes L^{-1}).
$$ 
and which implies that the generic curve $(C,B,L)$ in $\cR_{3,2}$ is a non-tetragonal curve. 

\end{proof}

\subsection{Proof of Theorem \ref{birational}}
Consider the inclusion $\cR_{3,2}\subset \cM_7$ of moduli spaces.
Then we recall the classification of the singular loci of the moduli space $\cM_g$ done by Cornalba \cite{Cornalba2}. In particular, 
the curves with non-trivial automorphisms lie in the singular locus of $\cM_g$ and precisely form the singular locus. The maximal components 
of the singular locus are also described by him.
We recall his result when $g=7$ and for the embedding $\cR_{3,2}\subset \cM_7$, since it will be crucial for us. We note that any double 
cover corresponding to $(C,B,L)\in \cR_{3,2}$ corresponds to an involution $i$ on $C'$ with four fixed points, and having the quotient $C=C'/i$.

\begin{proposition}\label{cornalba2}
The singular locus $\cS\subset \cM_7$ consists of smooth curves with automorphisms. 
In particular the moduli space $\cR_{3,2}$ lies in the singular locus $\cS$ and furthermore it is a maximal component of $\cS$.
\end{proposition}
\begin{proof}
See \cite[Corollary 1, p.146 and p.150]{Cornalba2}.
\end{proof}

Now, consider a generic point $(C'\sta{f}{\rar} C)=(C,B,\cL)\in \cR_{3,2}$. Then, by \cite[\S 2]{BCV}, we have a decomposition of the 
canonical space of $C'$:
\begin{equation}\label{eigenspaces}
H^0(C',\omega_{C'})= H^0(C,\omega_{C})\oplus H^0(C,\omega_{C}\otimes \cL).
\end{equation}
We can also interpret this decomposition for the involution $i$, which acts on the canonical space nontrivially.
Namely, we have a natural identification of the eigenspaces for $i$:
\begin{eqnarray*}
H^0(C',\omega_{C'})^+ & = & H^0(C,\omega_{C}) \\
H^0(C',\omega_{C'})^- & = & H^0(C,\omega_{C}\otimes \cL).
\end{eqnarray*}
Note that $\m{dim }H^0(C',\omega_{C'})^+=3$ and  $\m{dim }H^0(C',\omega_{C'})^-=4$.

We can now apply Theorem \ref{automorphism} to the automorphism $i$ and conclude that $i$ lifts to an automorphism $i$ of $\p^{15}$ and 
leaves $X_{10}$ invariant. This gives an action of $i$ on the representation space $U_{16}$. Indeed, since $\m{Pic}(X_{10})\simeq \Z$, 
the ample line bundle $\cO_{X_{10}}(1)$ is invariant under $i$. Hence $i$ induces an action on the sections of $ \cO_{X_{10}}(1)$ which is precisely $U_{16}$.
Let us write the eigenspace decomposition of $U_{16}$ for the $i$-action:
\begin{equation}\label{pmeigenspace}
U_{16}\,=\,U^+ \oplus U^-.
\end{equation}

There are various possibilities for the dimensions of $U^+$ and $U^-$, which will
correspond to
\begin{equation}\label{dimension}
(\m{dim }U^+,\m{dim }U^-):=(r,16-r), \m{ for } 1\leq r \leq 15,
\end{equation}
since $i$ acts nontrivially.

We make the following observation first.
\begin{lemma}\label{invariant}
A point of the product variety $G(3,U^+)\times G(4,U^-)\subset G(7,U_{16})$ corresponds to a linear space
$\p^6\subset \p^{15}$, which is $i$ invariant. Furthermore, if $\p^6$ intersects
$X_{10}$ transversely then the intersection is a non-tetragonal curve with an involution and satisfying the decomposition \eqref{eigenspaces}.
\end{lemma}
\begin{proof}
We first note that a $3$-dimensional subspace $V^+\subset U^+$ and $4$-dimensional subspace $V^-\subset U^-$ give a linear subspace 
$\p^6\subset \p^{15}$. Clearly $\p(V^+\oplus V^-)\subset \p(U_{16})$ is a $\p^6$ and is invariant under the action of $i$. For the 
second assertion, note that $C'=\p^6\cap X_{10}$ also is an $i$-invariant subset and whenever the intersection is transverse, it 
corresponds to a genus $7$ curve $C'$(by \cite{Mukai3}) with an involution, such that $\p^6$ is the canonical linear system of $C'$. 
This means that the $\pm$-eigenspaces of the canonical space of $C'$ are precisely $V^+$ and $V^-$. These data recover the decomposition 
in \eqref{eigenspaces}.
\end{proof}

\begin{lemma}\label{SO}
There is a subgroup $H\subset SO(10)$ such that $U^+$ and $U^-$ are $H$-representations. This induces an action of 
$H$ on $G(3,U^+)\times G(4,U^-)$ and which commutes with the action of $i$ such that the group quotient under this action is a
birational model of $\cR_{3,2}$. 
\end{lemma}
\begin{proof}
We note that by Mukai's classification \cite[\S 5]{Mukai3}, we have a birational isomorphism
$$
\cM_7  \sim  G(7,U_{16})/SO(10).
$$
The product subvariety $G(3,U^+)\times G(4,U^-)\subset G(7,U_{16})$ is  acted on not by $SO(10)$ but by an algebraic subgroup $H \subseteq SO(10)$. 
To describe the action of $H$, we first note that the involution $i$ commutes with the action of $H$, so that the quotient 
$(G(3,U^+)\times G(4,U^-))/H$ gives the isomorphism classes of smooth curves with an involution $i$.
Then the matrices in $SO(10)$ which act on the product subvariety are those which commute with the involution $i$ on a linear space 
$\p(U_{16})$.

As noted in \eqref{pmeigenspace}, we have an eigenspace decomposition
$$
U_{16}=U^+\oplus U^-
$$
for the action of $i$.
Since for any $h\in H$ and $s\in U^+$ (or $s\in U^-$)
$$
i.h(s)=h.i(s)=h(s)
$$
it follows that $U^+$ (or $U^-$) are (projective) $H$-modules (i.e, are projective representations of $H$ and hence $H$ acts on $\p(U^{\pm})$).

By Corollary \ref{generic}, we know that a generic curve $C'\in \cR_{3,2}$ is non-tetragonal. Hence, the moduli space
$\cR_{3,2}$ does not lie in the indeterminacy locus of the birational map
$$
\cM_7  \dasharrow  G(7,U_{16})/SO(10).
$$
Hence this birational map restricts to a generically injective rational map
$$
\psi: \cR_{3,2} {\dasharrow}  G(7,U_{16})/SO(10).
$$
Corresponding to a non-tetragonal curve $C'$, for $(C,B,L)=(C',i)\in \cR_{3,2}$ (which is the generic situation, by Corollary \ref{generic}) we can 
associate a point in $G(3,U^+)\times G(4,U^-)$ according to the decomposition of the canonical space in \eqref{eigenspaces}.
Hence the image of $\psi$ maps to the product space
$$
{\psi'}:\cR_{3,2}\sta{\psi'}{\dasharrow} (G(3,U^+)\times G(4,U^-))/H,
$$
and this map is generically injective.

To see that $\psi'$ is birational, given a generic point
$\p^6 \in G(3,U^+)\times G(4,U^-)$ we first know by \cite{Mukai3} that the intersection
$C'=\p^6\cap X_{10}$ lies in $\cM_7$. Now by Proposition \ref{cornalba2}, $C'$ lies in the singular locus $\cS\subset \cM_7$, since it 
has a nontrivial involution. This implies that the inverse image of  $(G(3,U^+)\times G(4,U^-))/H$ under $\psi$
in $\cM_7$ is a subset in the singular locus $\cS\subset \cM_7$ and containing a dense open subset of $\cR_{3,2}$.
But again by Proposition \ref{cornalba2} since $\cR_{3,2}$ is a maximal component in $\cS$, the inverse image has to be dense in $\cR_{3,2}$.

This proves the birational equivalence
\begin{equation}\label{something}
\cR_{3,2}\sim (G(3,U^+)\times G(4,U^-))/H.
\end{equation}
\end{proof}

\begin{corollary}
The moduli space $\cR_{3,2}$  (and hence ${\mathcal A}_4$) is a unirational variety.
\end{corollary}
\begin{proof}
Since a Grassmannian variety is a rational variety, it follows that the product space $G(3,U^+)\times G(4,U^-)$ is also
a rational variety. Using the description in \eqref{something}, it follows that the moduli space $\cR_{3,2}$ is a unirational
variety.
\end{proof}

The birational model in \eqref{something} should also be compatible with the projection $\cR_{3,2}\rar \cM_3$. 
It would be interesting to study $H$ in detail and therefore we pose the following question:

\begin{Question}\label{question}
 Determine the subgroup $H$ and the $H$-(projective) representations $U^+$ and $U^-$ explicitly.
\end{Question}

Notice that we have the spinor representation
$$
\phi(10):Spin(10)\rar Aut(U_{16})
$$
which gives the $SO(10)=\f{Spin(10)}{\pm 1}$-action on $P(U_{16})$, considered in \cite{Mukai3}.
It may be possible to study $H$ further via the spinor representation restricted to the various subgroups of $SO(10)$.

\subsection{Towards the motive of $\cR_{3,2}$ }

In this subsection we want to give one application concerning the motive of $\cR_{3,2}$. 

In \cite{Iy-Mu} we have constructed Chow--K\"unneth decompositions for open subsets of moduli space of curves of small genus $g\leq 8$.
Recall that this was proved in \cite{Iy-Mu}, via realizing the open subsets as group quotients of open subsets in homogeneous spaces.
The key point used was that the homogeneous spaces have only algebraic cohomology and hence orthogonal projectors equivariant for the 
group action could be constructed. All those methods can also be applied to the variety $\cR_{3,2}$. Using the birational equivalence in \eqref{something}, we obtain: 

\begin{proposition}
There is an open subset (not necessarily affine) of the moduli space $\cR_{3,2}$ which admits a Chow--K\"unneth decomposition.  
\end{proposition}

$\Box$



\begin{thebibliography}{AAAAA}
\bibitem[BCV]{BCV} F. Bardelli, C. Ciliberto, A. Verra, {\em Curves of minimal genus on a general abelian variety}, Compositio Math.  96  (1995),  no. \textbf{2}, 115--147.
\bibitem[B-dC]{B-dC} F. Bardelli, A. Del Centina, {\em Bielliptic curves of genus three: canonical models and moduli space}, 
Indag. Math. (N.S.) 10 (1999), no. \textbf{2}, 183--190.
\bibitem[Be]{Behrend} K. Behrend,  {\em On the de Rham cohomology of differential and algebraic stacks},  Adv. Math.  198  (2005),  no. \textbf{2}, 583--622.
\bibitem[Ch-Rn]{chen} M.C. Chang, Z. Ran, {\em Unirationality of the moduli spaces of curves of genus $11,$ $13$ (and $12$)}, 
Invent. Math. 76 (1984), no. \textbf{1}, 41--54.
\bibitem[Co]{Cornalba2} M. Cornalba, {\em On the locus of curves with automorphisms},  (Italian summary)
Ann. Mat. Pura Appl. (\textbf{4}) 149 (1987), 135--151.
\bibitem[Co2]{cornalba} M. Cornalba, {\em Cohomology of moduli spaces of stable curves}, 
Proceedings of the International Congress of Mathematicians, Vol. II (Berlin, 1998).  Doc. Math.  1998,  Extra Vol. II, 249--257.
\bibitem[Do]{Dolgachev} I. Dolgachev, {\em Rationality of $\cR\sb 2$ and $\cR\sb 3$}, Pure Appl. Math. Q. 4 (2008), no. \textbf{2}, part 1, 501--508.


\bibitem[Fa-Ve]{Fa-Ve} Farkas, G., Verra, A. {\em The classification of universal Jacobians over the moduli space of curves}, arXiv math.AG. 1005.5354.

\bibitem[Fu-Ha]{Fulton} W. Fulton, J. Harris, {\em Representation theory. A first course}, Graduate Texts in Mathematics, 129. Readings in Mathematics. Springer-Verlag, New York, 1991. xvi+551 pp.



\bibitem[Hn]{Arbarello} R. Hartshorne, {\em Algebraic geometry}, Graduate Texts in Mathematics, No. \textbf{52}. Springer-Verlag, New York-Heidelberg, 1977. xvi+496 pp.
\bibitem[Iy-Ml]{Iy-Mu} J. N. Iyer, S. M\"uller-Stach, {\em Chow--K\"unneth decomposition for some moduli spaces}, 
Documenta Mathematica, \textbf{14}, 2009, 1-18.
\bibitem[Iz-L-S]{Izadi} E. Izadi, M. Lo Giudice, G. Sankaran, {\em The moduli space of \'etale double covers of genus 5 curves is unirational}, Pacific J. Math. 239 (2009), no. \textbf{1}, 39--52.
\bibitem[Ka]{Katsylo} P. Katsylo, {\em Rationality of the moduli variety of curves of genus $3$}, Comment. Math. Helv. 71 (1996), no. \textbf{4}, 507--524.
\bibitem[Mk]{Mukai3} S. Mukai, {\em Curves and symmetric spaces. I}, Amer. J. Math.  117  (1995),  no. \textbf{6}, 1627--1644.
\bibitem[Mk2]{Mukai6} S. Mukai, {\em Curves and symmetric spaces},  Proc. Japan Acad. Ser. A Math. Sci.  68  (1992),  no. \textbf{1}, 7--10.
\bibitem[Ra]{Ramanan} S. Ramanan, {\em Ample divisors on abelian surfaces}, 
Proc. London Math. Soc. (3)  51  (1985),  no. \textbf{2}, 231--245.

\bibitem[Se]{Sernesi} E. Sernesi, {\em Unirationality of the variety of moduli of curves of genus twelve}, 
(Italian) Ann. Scuola Norm. Sup. Pisa Cl. Sci. (4) 8 (1981), no. \textbf{3}, 405--439.


\bibitem[Ve]{Verra} A. Verra, {\em The unirationality of the moduli spaces of curves of genus 14 or lower}, 
Compos. Math. 141 (2005), no. \textbf{6}, 1425--1444.

\bibitem[Ve2]{Verra2} A. Verra, {\em On the universal principally polarised abelian variety of dimension 4}, 
Curves and abelian varieties, 253--274, Contemp. Math., \textbf{465}, Amer. Math. Soc., Providence, RI, 2008. 

\end{thebibliography}
\end{document}